\documentclass{amsdtx}
\usepackage{amsmath}        
\usepackage{amsfonts}
\usepackage{amscd}
\usepackage{amssymb}

\usepackage{graphicx}
\usepackage{xypic}
\input xypic
\input xy
\xyoption{all}

\newcommand{\lga}{\longrightarrow}
\newcommand{\la}{\langle}
\newcommand{\ra}{\rangle}

\begin{document}
\title{A note on the characteristic classes of non-positively curved manifolds}
\author{Jean-Fran\c{c}ois Lafont, Ranja Roy}
\date{ }
\maketitle
\centerline{\it Department of Mathematical Sciences, Binghamton University, Binghamton, NY-60902}
\centerline{\it email : jlafont@math.binghamton.edu}
\centerline{\it Department of Mathematics, New York Institute of Technology,
Old Westbury, NY-11568}
\centerline{\it email : rroy@nyit.edu}
\noindent
------------------------------------------------------------------------------------------------------------------------------------
\\
{\bf Abstract }\\
\medskip
 
\noindent 
We give conceptual proofs of some well known results concerning compact non-positively curved locally symmetric spaces. 
We discuss vanishing and non-vanishing of Pontrjagin numbers and Euler characteristics
for these locally symmetric spaces.  We also establish vanishing results for Stiefel-Whitney numbers of (finite covers
of) the
Gromov-Thurston examples of negatively curved manifolds.  We mention some geometric corollaries: the MinVol
question, a lower bound for degrees of covers having tangential maps to the non-negatively curved duals,
and estimates for the complexity of some representations of certain uniform lattices.
\medskip

\medskip

\noindent
{\sl AMS primary classification:} 57R20, 53C35.
{\sl AMS secondary classification:} 57T10.\\
\noindent
------------------------------------------------------------------------------------------------------------------------------------

\vspace{0.2 in}

\noindent
{\bf 1. Introduction.}
\vskip 10pt

A well known result asserts that closed hyperbolic manifolds have zero Pontrjagin numbers. The 
standard argument for this consists of observing that closed hyperbolic manifolds are conformally 
flat, and that conformally flat closed manifolds have zero Pontrjagin numbers by Avez \cite{av}. 
This note originated in a desire to give a simple, conceptual proof of 
this basic result, which we do in section 2.  Our argument gives an alternate, more geometrical viewpoint
that should be contrasted with the classical
Hirzebruch proportionality principle.  The main advantage of our approach lies in that the characteristic numbers are
computed {\it via an actual map} between the non-positively curved locally symmetric spaces and their non-negatively
curved duals. 

Now recall that there is another well-known class of negatively curved closed manifolds arising from 
the Gromov-Thurston construction \cite{gt} (see also the older construction of Mostow-Sui \cite{ms}). 
These manifolds are ramified coverings of closed hyperbolic 
manifolds, where the ramification occurs over a totally-geodesic, codimension two submanifold that 
is null-homologous. Note that the behavior of characteristic numbers under ramified coverings 
is unclear (though see the recent result of Izawa \cite{iz}). In section 3, we show that the Gromov-Thurston 
manifolds always have a finite cover that bounds orientably. A 
byproduct of our argument also gives a very elementary proof of a result of Rohlin \cite{ro}: 
Every orientably closed 3-manifold bounds orientably.

Finally, in section 4, we point out various geometric corollaries of our main results. While many of these are 
standard, we do include some new results.  We conclude in section 5 with some open questions.

\vskip 10pt

\noindent
{\bf 2. Characteristic numbers of negatively curved locally symmetric spaces. }
\vskip 10pt

Let us start by recalling the construction of the non-negatively
curved dual space associated to any non-positively curved closed locally symmetric space. 
If $G$ is a real Lie group, $K$ it's maximal compact subgroup, we let 
$G_C = G \otimes {\mathbb C}$ be the complexification of $G$ and $G_{U}$ the maximal compact subgroup 
of $G_C$. The factor spaces $G/K$ and $M_U = G_U/K$ are called dual symmetric spaces \cite{bo}.  By abuse
of language, if $\Gamma$ is a uniform lattice in $G$, we will still say that $M:= \Gamma \backslash G/K$
and $M_U$ are dual spaces.  In \cite{bo}, Okun showed that if $M^n$ is a closed locally symmetric space, 
then there is a tangential map from some finite cover $\bar M^n$ to the dual symmetric space.  We start
by showing the following easy Lemma:

\vskip 5pt

\noindent
{\bf Lemma 1:} Assume $f: M\longrightarrow N$ is a tangential map between two n-dimensional manifolds. Then
\begin{itemize}
\item $p_I(M) = \pm \mbox{deg}(f) \cdot p_I(N) \in {\mathbb Z}$
\item $sw_I(M) = \mbox{deg}(f) \cdot sw_I(N) \in {\mathbb Z}/ 2{\mathbb Z}$
\end{itemize}
where $p_I$, $sw_I$ denote the Pontrjagin and Stiefel-Whitney numbers associated to a product of Pontrjagin
or Stiefel-Whitney classes.

\vskip 5pt

\noindent
{\bf Proof:} Since the map is tangential, the pullbacks of Pontrjagin classes (respectively Stiefel-Whitney 
classes) of $N$ yield the corresponding classes for $M$.
If we denote by $\tau _I(N)$ a product of Pontrjagin classes, we have $f^{*}(\tau _I(N)) = \tau _I(M)$. 
Likewise, if $\sigma _I(N)$ denotes a product of Stiefel-Whitney
classes, $f^{*}(\sigma _I(N)) = \sigma _I(M)$. Now we have that:\\
\centerline{$p_I(M) = \la \tau _I(M), [M] \ra =  \la f^{*}(\tau _I(N)), [M] \ra $}
\centerline{$=\pm \la \tau _I(N), f_{*}([M]) \ra =  \pm \la (\tau _I(N)), \mbox{deg}(f)\cdot [N] \ra $}
\centerline{$=\pm \mbox{deg}(f)\cdot \la (\tau _I(N)), [N] \ra = \pm \mbox{deg}(f)\cdot p_I(N)$.}
\noindent
And the argument for part (b) of the lemma is identical.

\vskip 5pt

Note that, from the discussion above, we have associated to any closed locally symmetric space $M^n$ a 
diagram:
$$M^n \longleftarrow \bar{M}^n \longrightarrow M_U$$
where $\bar{M} ^n$ is a finite cover, $M_U$ is the non-negatively curved dual, and the maps in the diagram
are tangential.  Since a covering map never has zero degree, Lemma 1 tells us that we can solve for
the Pontrjagin numbers of $M^n$:

$$p_I(M^n)=\frac{deg(t)}{deg(f)}\cdot p_I(M_U)$$

Of course, if we are trying to obtain vanishing/non-vanishing of Pontrjagin numbers, it is crucial to know
when $deg(t)\neq 0$.  Conceivably if $deg(t)=0$, one could have non-zero Pontrjagin numbers for $M_U$, but with
the corresponding Pontrjagin number for $M^n$ equal to zero.  That this does not occur is the content of the 
next Lemma:

\vskip 5pt

\noindent
{\bf Lemma 2:} If $t$ has degree zero, then the Pontrjagin numbers $p_I(M_U)$ are all equal to zero.
\vskip 5pt

\noindent
{\bf Proof:} We start by noting that Okun (\cite{bo}, Corollary 6.5) showed that if $G_U$ and $K$ have equal
rank, then $t$ has non-zero degree.   Hence if $deg(t)=0$, we must have $rk(G_U)>rk(K)$.  Recall that
the {\it toral rank} of a compact manifold $N$, denoted by $trk(N)$, is the largest dimension of a torus that 
has a smooth, 
rationally-free action on $N$ (see \cite{sp}). Now Allday-Halperin \cite{ah} have shown that $trk(G_U/K)=
rk(G_U)-rk(K)$, hence if $deg(t)=0$, we have that $trk(M_U)>0$.  But Conner-Raymond \cite{cr} have shown that
if $N$ is a compact manifold with $trk(N)>0$, then all the Pontrjagin numbers of $N$ are equal to zero.  Applying
their result to $M_U$ completes the proof.

\vskip 5pt

For completeness, we point out that by a result of Papadima \cite{sp}, for the homogenous space $M_U=G_U/K$,
we have that the toral rank of $M_U$ is zero if and only if the Euler characteristic of $M_U$ is non-zero.
Hence to verify that the map $t$ has non-zero degree, it is sufficient to verify that the Euler characteristic
of $M_U$ is non-zero.  
We refer to Helgason \cite{hl} for the classification of the irreducible higher rank non-positively curved
closed locally symmetric spaces, as well as for the notation used in our discussion.  
The results for the classical families can be summarized in the following:

\vskip 5pt

\noindent 
{\bf Theorem A}: Let $M^n$ be a closed orientable irreducible higher rank locally symmetric space, and assume that $M^n$
is locally modeled on one of the following:
\begin{enumerate}
\item $SU(p,q)/S(U_p\times U_q)$
\item $SO_0(p,q)/SO(p)\times SO(q)$ where $p$ and $q$ are not both odd
\item $SO^*(2n)/U(n)$
\item $Sp(n, \mathbb{R})/U(n)$
\item $Sp(p,q)/S(p)\times S(q)$
\end{enumerate}
Then $M^n$ has non-zero Euler characteristic.  
\vskip 5pt

\noindent
{\bf Theorem B}: Let $M^n$ be a closed orientable irreducible higher rank locally symmetric space, and assume that $M^n$
is locally modelled on one of the following:
\begin{enumerate}
\item $SL(n, \mathbb {R})/SO(n)$
\item $SU^*(2n)/Sp(n)$
\item $SO_0(p,q)/SO(p)\times SO(q)$ where $p$ and $q$ are both odd
\item an irreducible globally symmetric spaces of Type IV, see pgs. 515-516 in \cite{hl}
\end{enumerate}
Then $M^n$ has all Pontrjagin numbers equal to zero.

\vskip 5pt

\noindent
{\bf Proof of Theorems A \& B:}  Let us first note that the duals of the locally symmetric spaces mentioned in 
theorem A are respectively:
\begin{enumerate}
\item $SU(p+q)/S(U_p\times U_q)$
\item $SO(p+q)/SO(p)\times SO(q)$ where $p$ and $q$ are not both odd
\item $SO(2n)/U(n)$
\item $Sp(n)/U(n)$
\item $Sp(p+q)/S(p)\times S(q)$
\end{enumerate}
Since the ranks of the various Lie groups being considered above are $rk(Sp(n))=n$, $rk(SO(2n))=rk(SO(2n+1))=n$,
$rk(SU(n))=n-1$, $rk(U(n))=n$, and $rk(S(U_p\times U_q))=p+q-1$, we see that in all the cases of Theorem A, the
condition $rk(G_u)=rk(K)$ is satisfied.  Since the toral rank of $M_u$ is zero, this implies that the Euler characteristic
of $M_U$ is non-zero, giving Theorem A.

Likewise for the cases appearing in Theorem B, we have that the duals are respectively:
\begin{enumerate}
\item $SU(n)/SO(n)$
\item $SU(2n)/Sp(n)$
\item $SO(p+q)/SO(p)\times SO(q)$ where $p$ and $q$ are both odd
\item a Lie group
\end{enumerate}
In each of the first three cases, we see that $rk(G_U)>rk(K)$.  So by the argument in Lemma 2, we have that all the 
Pontrjagin numbers of the dual spaces $M_U$ are zero.  Hence the Pontrjagin numbers for $M^n$ are likewise zero.  For the
fourth case, we note that Lie groups are parallelizable, hence have all Pontrjagin numbers zero. This concludes the proof.

\vskip 5pt

\noindent
{\it Remark:}  Theorem A above lists the only irreducible higher hank locally symmetric spaces of non-positive 
curvature which could conceivable have non-vanishing Pontrjagin numbers.  Since the procedure for calculating the
Pontrjagin numbers of the non-negatively curved duals is well established (see Borel-Hirzebruch \cite{bh}), one could in 
principle find out which of these spaces actually have a non-vanishing Pontrjagin number (note that by Lemma 2, for the
spaces in Theorem A, the degree of the tangential map is non-zero).  As this procedure is primarily 
combinatorial in nature, we leave the precise computations to the interested reader, and content ourselves with
computing them for the {\it negatively} curved locally symmetric spaces.  In the process, we also discuss the 
exceptional
locally symmetric space $F_{4(-20)}/Spin(9)$ giving rise to Cayley hyperbolic manifolds.  We leave to the interested
reader the task of deciding for the remaining exceptional cases which of Theorem A or B applies.

\vskip 5pt

\noindent 
{\bf Corollary 1}: Let $M^n$ be a compact orientable manifold, and assume that one of the following holds:
\begin{enumerate}
\item $M^n$ is real hyperbolic
\item $M^n$ is complex hyperbolic, and $n=4k+2$
\item $M^n$ is quaternionic hyperbolic, and $n=8k+4$
\end{enumerate}
Then $M^n$ has a finite cover that bounds.  In the first two cases, there is a finite cover that bounds
{\it orientably} (and hence $M^n$ has all Pontrjagin numbers equal to zero).
\vskip 5pt

\noindent 
{\bf Corollary 2}: Let $M^n$ be a compact orientable manifold, and assume that one of the following holds:
\begin{enumerate}
\item $M^n$ is Cayley hyperbolic (so $n=16$)
\item $M^n$ is complex hyperbolic, and $n=4k$
\item $M^n$ is quaternionic hyperbolic of dimension at least $8$.
\end{enumerate}
Then $M^n$ has some non-zero Pontrjagin numbers, and hence no finite cover can bound orientably.  Furthermore,
in the first two cases, we have that {\it all} Pontrjagin numbers are non-zero.

\vskip 5pt

\noindent
Since the arguments are closely related, we simultaneously prove both corollaries:

\vskip 5pt

\noindent
{\bf Proof of Corollaries 1 \& 2}: We note that for the negatively curved symmetric spaces, the duals are easy 
to compute.  Indeed we have that:
\begin{itemize}
\item the dual to real hyperbolic space is the sphere,
\item the dual to complex hyperbolic space is complex projective space,
\item the dual to quaternionic hyperbolic space is quaternionic projective space,
\item the dual to Cayley hyperbolic space is the Cayley projective plane.
\end{itemize}
Since the characteristic classes
of the duals are well known, we can apply Lemmas 1 and 2 in each case to obtain information on the negatively
curved locally symmetric spaces.  $\bar M^n$ will always denote the finite cover that supports a tangential
map to the positively curved dual.  The various cases are:

\vskip 5pt

\noindent
{\bf $M^n$ is real hyperbolic}: Since the sphere bounds orientably, all its characteristic numbers (both Stiefel-Whitney
and Pontrjagin) are zero.  Applying Lemma 1,
we see that all the characteristic numbers of $\bar M^n$ are zero. By a result of Wall \cite{wl}, this is equivalent
to $\bar M^n$ bounding orientably, giving (1) of Corollary 1.

\vskip 5pt

\noindent
{\bf $M^{2n}$ is complex hyperbolic}: Then its dual space is the complex projective space ${\mathbb C}P^n$, which is 
a 2n-dimensional real manifold.  We now have two cases:

\noindent
(A) If $n= 2k$, then the Pontrjagin numbers are all non-zero (\cite{mi} page 185), hence using Lemmas 1 and 2, 
the same holds for $M^{2n}$.

\noindent
(B) If $n=2k+1$, then ${\mathbb C}P^n$ bounds orientably (\cite{mi} page 186). Arguing as in the real hyperbolic 
case, we see that $\bar M^n$ bounds orientably.

\noindent This gives us (2) of Corollaries 1 and 2.

\vskip 5pt

\noindent
{\bf $M^{4n}$ is quaternionic hyperbolic}: Then its dual space is the quaternionic projective space ${\mathbb O}P^n$, which is 
a 4n-dimensional real manifold.  We again have two cases:

\noindent
(A) If $n=2k+1$, then ${\mathbb O}P^n$ bounds, and hence has vanishing Stiefel-Whitney numbers.  By
Lemma 1, the same holds for $\bar M^{2n}$, giving (3) of Corollary 1.


\noindent
(B) In general, the total Pontrjagin class of $\mathbb OP^n$ is given by $(1 + u)^{2n+2}(1 + 4u)^{-1}$, where 
$u\in H^4({\mathbb O}P^n)$
is a generator for the truncated polynomial ring $H^*({\mathbb O}P^n)$.  Since the coefficient of $u$ in the power series 
expansion equals $2n-2$, we see that the Pontrjagin number $p_1^{n}(M_U)$ is
equal to $(2n-2)^{n}$.  So provided $n\geq 2$, we can apply Lemmas 1 and 2 to obtain (3) of Corollary 2.

\vskip 5pt

\noindent
{\bf $M^{16}$ is Cayley hyperbolic}: Then its dual space is the Cayley projective plane 
$\mbox{Cay}P^2$. The Cayley plane has two non-vanishing Pontrjagin numbers, namely 
${p_2}^2[\mbox{Cay}P^2] = 36$ and ${p_4}[\mbox{Cay}P^2] = 39$ (see Borel-Hirzebruch \cite{bh}, pages 535-536). 
Applying Lemma 2, we get that $\bar M^{16}$ has some non-vanishing Pontrjagin numbers.  This deals with case (1)
of Corollary 2, and hence completes the proof of the Corollaries.

\vskip 5pt

\noindent 
{\it Remark}: We note that information on the Stiefel-Whitney numbers of the rank one locally symmetric spaces is
much harder to obtain.  Indeed, anytime the degree of one of the two maps is even, there is a potential loss of 
information.  

\vskip 5pt

\noindent
{\bf Corollary 3:}  If $M^n$ is a manifold supporting a metric of constant sectional curvature, then all of its 
Pontrjagin numbers are zero.

\vskip 5pt

\noindent
{\bf Proof}:  The case of constant negative curvature has been dealt with above.  In the remaining two cases, $M^n$
has a finite cover that bounds orientably (either a sphere, or a torus, depending on curvature).  The corollary follows.

\vskip 5pt

\noindent 
{\it Remark}: Recall that Farrell-Jones have constructed exotic smooth 
structures on certain closed hyperbolic manifolds, and have shown that 
these manifolds support Riemannian metrics of negative curvature \cite{fj}.  There results were subsequently extended to
providing exotic smooth structures on a variety of different locally symmetric spaces, see for instance \cite{fj2}, \cite{po},
\cite{fjo}, \cite{af2}, \cite{af}.  Observe that while the Pontrjagin classes are
smooth invariants, the rational Pontrjagin classes are topological invariants, by a celebrated result of Novikov \cite{nov}.
Since the Pontrjagin numbers of a manifold only depend on the rational Pontrjagin classes (i.e. the torsion part of the
Pontrjagin classes do not influence the Pontrjagin numbers), the discussion in Theorems A \& B gives us vanishing 
(or non-vanishing) results for the Pontrjagin numbers of these exotic manifolds as well.  

\vskip 10pt

\noindent
{\bf 3. Characteristic numbers of the Gromov-Thurston examples.  }
\vskip 10pt

\noindent
{\bf Definition:} Let $X$ be a oriented differentiable manifold (with or without boundary) on which the cyclic 
group $\mathbb{Z}_k$
acts semifreely by orientation-preserving diffeomorphisms with fixed set a codimension two submanifold $Y$.
Denote the quotient space by $X^\prime:=X/\mathbb{Z}_k$, and the canonical projection map by $\pi:X\rightarrow 
X^\prime$.  Let $Y$ be the fixed set of the action on $X$, and note that $\pi:Y\rightarrow Y^\prime$ is a diffeomorphism.
Observe that $X^\prime$ is a manifold.  We say that the $X$ is an {\it oriented cyclic ramified cover} of $X^\prime$, of
order $k$, ramified over $Y^\prime$.  If $Y^\prime$ bounds a smooth embedded codimension one submanifold in $X^\prime$, 
we say that the ramified covering is {\it nice}. 

\vskip 5pt
%
%

\noindent
{\it Remark:} If a ramified covering is nice, then it is particularly easy to describe it.  Indeed, let $N$ be the
codimension one embedded submanifold satisfying $\partial N=Y^\prime$.  Then the pre-image of $N$ in the ramified
cover $X$ will consist of multiple (embedded) copies of $N$ which all coincide along their boundary (which will 
equal $Y$).  Cutting $X$ open along the pre-images of $N$ will yield $k$ homeomorphic copies of $X^\prime -N$.  Now 
consider the space with boundary  the double $DN$ of $N$, obtained by cutting open $X^\prime$ along $N$.  Then $X$
is obtained by taking $k$ copies of this space, $X_1, \ldots X_k$, and for each space, cyclically gluing $\partial X_i^+$
to $\partial X_{i+1}^-$, where $\partial X_i^\pm$ denotes the two copies of $N$ in $\partial X_i=DN$.

\vskip 5pt

\noindent
{\bf Proposition:} Assume that $M^n$ bounds, and that $p:\bar M^n\rightarrow M^n$ is an oriented
cyclic ramified cover of $M^n$ (ramified over $N^{n-2}$).  If the covering is nice, then $\bar M^n$ also bounds.
If $M^n$ bounds orientably, then so does $\bar M^n$.

\vskip 5pt

\noindent
{\bf Proof:} Let $M^n = \partial L^{n+1}$, and note that since the ramified covering is nice, there exists a 
a smoothly embedded $K_0^{n-1}\subset M^n$ satisfying $\partial K_0^{n-1}=N^{n-2}$. Since $M^n$ is collarable in 
$L^{n+1}$, there is a manifold $K^{n-1}\subseteq L^{n+1}$ of dimension ${n-1}$ with the properties:
\noindent
\begin{itemize}
\item $K^{n-1}\cap \partial L^{n+1} = N^{n-2} = \partial K_0^{n-1}$,
\item $K^{n-1}$ and $K_0^{n-1}$ are cobordant in $L^{n+1}$,
\item the cobordism $W^n$ is an embedded submanifold satisfying $W^n\cap M^n = K_0^{n-1}$.
\end{itemize}
\noindent
Indeed, homotoping $K_0^{n-1}$ (relative $\partial K_0^{n-1}=N^{n-2}$) into a collared neighborhood of $M^n$ in $L^{n+1}$
give both $K^{n-1}$, and the manifold $W^n$ (the image of the homotopy, which we can assume to have no self-intersections).
Now note that $K^{n-1}\subseteq L^{n+1}$ is a codimension two submanifold which bounds $W^n$. Hence we can take the 
$i$-ramified covering of $L^{n+1}$ over $K^{n-1}$ (see the remark preceding this Proposition). But note that on 
$\partial L^{n+1}=M^n$, this ramified covering 
yields ${\bar M}^n$. Hence if ${\bar L}^{n+1}$ is the covering, we have $\partial{\bar {L}}^{n+1}= {\bar M}^n$.
Finally, we note that if $L^{n+1}$ is orientable, then so is the ramified covering $\bar L^{n+1}$.

\vskip 5pt

\noindent
{\bf Corollary 4:} Let $M$ be a closed, orientable, 3-dimensional manifold. Then $M$ bounds orientably. (This result is 
originally due to Rohlin \cite{ro} ).

\vskip 5pt

\noindent
{\bf Proof: } It is a well known result (due independantly to Hilden \cite{hi} and Montesinos \cite{mo}) that every 
closed orientable 3-manifold is ramified covering of the 3-dimensional 
sphere $S^3$ along a knot. Since every knot in $S^3$ bounds a compact embedded surface, this ramified cover is 
nice.  Since $S^3$ bounds orientably, the proposition gives us the claim.

\vskip 5pt

\noindent
{\it Remark:} The Corollary also follows easily from results of Thom and Wall: the Pontrjagin numbers are automatically
zero, since $M$ is 3-dimensional.  As for the Stiefel-Whitney numbers, there are only three of them to consider:
$s_1^3, s_1^2s_2$, and $s_3$.  Note that since $M$ is orientable, $s_1=0$, so the first two numbers vanish.  As for
$s_3$, it is just the mod $2$ reduction of the Euler characteristic, which has to be zero as we are in odd dimension.
Applying Wall's theorem (\cite{wl}), we get that $M$ must bound orientably.  The advantage of our approach is that
the bounding manifold can be seen {\it explicitly}, and we avoid appealing to the sophisticated results of Thom and 
Wall.

\vskip 5pt

\noindent
{\bf Theorem C:} Let $N$ be the Gromov-Thurston non-positively curved manifold. Then $N$ has a finite cover that 
bounds orientably (and hence all Pontrjagin numbers of $N$ are zero).

\vskip 5pt

\noindent
{\bf Proof of Theorem C:} Let $M$ be a real hyperbolic manifold and $N$ be the Gromov-Thurston non-positively curved manifold obtained
as a ramified covering of $M$. From Corollary 1, $M$ has a finite cover $\bar M$ that bounds orientably. 
We claim that there is a space $\bar N$ yielding the commutative diagram:
$$\xymatrix{\bar N \ar[r]^{\bar \psi} \ar[d]_{\bar \phi} & \bar M \ar[d]^{\phi}\\
N \ar[r]_{\psi} & M 
}$$
where $\bar \phi$ is a covering map and $\bar \psi$ is a 
ramified covering (and $\psi$ is the original ramified covering, $\phi$ the 
original covering).

In order to see this, we make the following general observation: assume that $X^{n-2}$ is a smooth embedded codimension 
two submanifold in $Y^n$, and let $W\subset Y^n$ be a closed tubular neighborhood of $X^{n-2}$.  Note that 
$W$ is a $\mathbb D^2$-bundle over $X^{n-2}$, and hence that $\partial W$ is an $S^1$-bundle over $X^{n-2}$.
Now let $Y^\prime \subset Y^n$ be the manifold with boundary obtained by removing the interior of $W$ from $Y^n$,
and assume that $\bar Y^\prime \rightarrow Y^\prime$ is a covering map.  Then we have:
\begin{enumerate}
\item the covering map $f:\bar Y^\prime \rightarrow Y^\prime$ extends to a covering $\bar f:\bar Y\rightarrow Y$ if and 
only if, for each fiber $F$ of the bundle $S^1\rightarrow \partial W\rightarrow X^{n-2}$, we have that $f^{-1}(F)$
consists of $deg(f)$ disjoint copies of $S^1$.
\item the covering map $f:\bar Y^\prime \rightarrow Y^\prime$ extends to a {\it ramified} covering 
$\bar f:\bar Y\rightarrow Y$ of degree $deg(f)$ over $X^{n-2}$ if and 
only if, for each fiber $F$ of the bundle $S^1\rightarrow \partial W\rightarrow X^{n-2}$, we have that $f^{-1}(F)$
is connected.
\end{enumerate}
Indeed, one direction of the implications is immediate, since a covering (respectively a ramified covering over 
$X^{n-2}$) exhibits precisely the aforementioned behavior on the boundary of a regular neighborhood.  Conversely,
assume that we have a covering map $f:\bar Y^\prime \rightarrow Y^\prime$ satisfying one of the above properties.
Then note that the pre-image $f^{-1}(\partial W)$ naturally inherits a smooth foliation with $S^1$ leaves.  Now 
consider the space $\bar W$ obtained by smoothly gluing in $\mathbb D^2$'s along their boundary to the leaves.  Observe
that this can be done, since the foliation on $f^{-1}(\partial W)$ is the lift of a fibration, and hence is locally 
a product.  Finally, form the space $\bar Y$ by gluing $\bar Y^\prime$ with $\bar W$ along their common boundary.

Now in case (1) above, we immediately get that the covering map $f$ extends to a covering map $\bar f$, by simply
extending linearly along each $\mathbb D^2$.  In case (2), we again extend linearly, but this time also 
extend the action of $\mathbb Z_{deg(f)}$ (by deck transformations) from each $S^1$ to each $\mathbb D^2$.  Note that
this gives a smooth $\mathbb Z_{deg(f)}$ action on $\bar Y$, whose fixed point set maps diffeomorphically to 
the original $X^{n-2}$.  

Now in the setting we have, proceed as follows: if $K^{n-2}$ is the codimension two submanifold of $M^n$ that is 
being ramified over, then let $W$ be a tubular neighborhood of $K$, $W_0$ it's interior.  Note that $\psi$ is 
an actual {\it covering}, when restricted to the preimage of $M-W_0$ (as we are throwing away a neighborhood of 
the set where the ramification occurs).  Consider the commutative 
diagram:
$$\xymatrix{M^\prime \ar[rr] \ar[d] & & \phi^{-1}(M-W_0) \ar[d]^{\phi}\\
\psi^{-1}(M-W_0) \ar[rr]_{\psi} & & M-W_0 
}$$
where $M^\prime$ is the pullback of the covering maps.  By commutativity of the diagram, we see that the covering
$M^\prime\rightarrow \phi^{-1}(M-W_0)$ satisfies (2) from our discussion above, while the covering 
$M^\prime\rightarrow \psi^{-1}(M-W_0)$ satisfies (1) from the discussion above.  In particular, extending $M^\prime$
as above, we obtain a space $\bar N$ which is simultaneously a ramified covering of $\bar M$, and an actual covering of 
$N$, as desired.

Finally, we note that the ramified covering $\bar \psi:\bar N\rightarrow \bar M$ is nice.  Indeed, in the Gromov-Thurston
construction, the ramified covering $\psi:N\rightarrow M$ is nice, so we have that $K^{n-2}=\partial L^{n-1}$ for a s
mooth, embedded codimension one manifold with boundary.
But we have that the map $\bar \psi$ is ramified over $\phi^{-1}(K^{n-2})$, which clearly bounds the smooth, embedded
codimension one submanifold $\phi^{-1}(L^{n-1})$.  This confirms that $\bar \psi$ is nice, and since $\bar M$ bounds
orientably, applying the Proposition, we see that $\bar N$ bounds orientably as well.  This completes the proof of 
Theorem C.

\vskip 5pt

\noindent
{\it Remark:} A related (unpublished) result is due to Ardanza-Trevijano Moras \cite{at}, and asserts 
that for the Gromov-Thurston
ramified coverings, the individual Pontrjagin classes vanish.  We note that while our
approach does not give vanishing of individual {\it classes}, it does give vanishing of the 
Stiefel-Whitney numbers
on a finite cover (which does not follow from the approach in \cite{at}).

\vskip 10pt

\noindent
{\bf 4. Geometric applications.}
\vskip 10pt

As is well known, characteristic numbers provide obstructions to a wide range of topological problems.  To mention
but a few, if $M^n$ has a non-zero Pontrjagin number, then 
\begin{enumerate}
\item no finite cover of $M^n$ bounds orientably.
\item $M^n$ has no orientation reversing self-diffeomorphism.
\item $M^n$ does not support an almost quaternionic structure (\cite{so}). 
\end{enumerate}
From our Corollary 2, we immediately get these properties for the rank one locally symmetric manifolds that are
either complex hyperbolic (with $n=4k$), quaternionic hyperbolic or Cayley hyperbolic.
\vskip 5pt

For another application, we note that while our Theorem A does not tell us which of the irreducible, higher-rank,
non-positively curved compact manifolds have non-zero Pontrjagin numbers, it does tell us which of these have non-zero
Euler characteristic.  Again, the Euler characteristic is known to be an obstruction to various topological/geometrical
problems, for instance for the spaces discussed in Theorem A, we immediately get that there cannot exist a nowhere
vanishing vector field.  For a possibly more interesting application, recall that the $MinVol$ of a smooth manifold is
defined to be the infimum of the volumes $Vol(M^n,g)$ over all Riemmanian metrics $g$ whose curvatures are bounded 
between $-1$ and $1$.  It is known that if $MinVol(M^n)=0$, then the Euler characteristic is also zero (see for instance
the survey \cite{pa}).  
This gives us:

\vskip 5pt

\noindent
{\bf Corollary 5:} Let $M^n$ be a compact orientable manifold which is locally symmetric.  If $M^n$ is
one of the manifolds allowed in Theorem A, then $MinVol(M^n)>0$.  

\vskip 5pt

Again, this result is known to hold for {\it all} non-positively curved compact locally symmetric spaces, and follows 
from the
estimates of filling invariants found in Gromov's work \cite{gr}.  For a more recent (and more general) proof 
using barycenter map techniques, see Connell-Farb \cite{cf}.

\vskip 5pt

\noindent
{\bf Corollary 6:} Let $M^{4n}$ be a compact orientable manifold which is locally symmetric.  Assume that $M^{4n}$ is
one of the manifolds allowed in Theorem A. For each partition 
$I = i_1,i_2,...,i_r$ of $n$, let $p_I(M^{4n})$ (respectively 
$p_I(M_U)$) denote the I-th Pontrjagin number of $M^{4n}$ (respectively of the dual $M_U$). Note 
that if $p_I(M_U) \neq 0$, then we also have that $p_I(M^{4n})\neq 0$ (from Lemma 2). Define 
$$\mu(M^{4n}) = LCM_I\{LCM(p_I(M^{4n}),p_I(M_U))/p_I(M^{4n})\}$$ 
where $LCM$ 
denotes least common multiple, and the outer $LCM$ is over all partitions $I$ of $n$ for which $p_I(M^{4n})\neq 0$. 
If $\bar{M}^{4n}\lga M^{4n}$ is 
a degree $d$ cover having a tangential map $\bar{M}^{4n}\lga M_U$, then $\mu(M^{4n})$ divides $d$.

\vskip 5pt

\noindent
{\bf Proof:} Let $r$ be the degree of the tangential map $\bar{M}^{4n}\lga M_U$. Then for each $I$, 
we have that $d\cdot p_I(M^{4n}) = r\cdot p_I(M_U)$. This implies that $d\cdot p_I(M^{4n})$ is a multiple of 
$LCM(p_I(M^{4n}),p_I(M_U))$. Hence for each $I$,we see hat $d$ is a multiple of 
$\frac{LCM(p_I(M^{4n}),p_I(M_U))}{p_I(M^{4n})}$. This forces $d$ to be a multiple of their least common 
multiple. Therefore $d$ is a multiple of $\mu(M^{4n})$.

\vskip 5pt

\noindent
{\it Remark:} The argument for the last corollary applies equally well to give an identical estimate for the
degree of the tangential map from $\bar M^n$ to $M_U$.  Part of our interest in the covering map (rather than
the tangential map), stems from the following:

\vskip 5pt

\noindent
{\bf Corollary 7:} Let $G/K$ be one of the irreducible globally symmetric spaces allowed as local models in the
hypotheses of Theorem A, and assume the dimension of $G/K$ is divisible by $4$.  Let $\Gamma $ be a torsion free subgroup 
of $G$, and denote by
$\Gamma \backslash G /K =: M^{4n}$ the associated locally symmetric space. Consider the flat principal bundle 
$G/K \times _{\Gamma } G \lga M^{4n}$, and extend its structure group to the group $G_C$. 
The bundle naturally defines a homomorphism $\rho : \Gamma \lga G_C\subset GL(k,\mathbb {C})$ (for some
suitable $k$). Let $A \subseteq \mathbb{C}$ 
be {\it any} subring of $\mathbb C$, finitely generated, with the property that 
$\rho (\Gamma )\subseteq GL(k, A)$, and let $m_1$, $m_2$ be {\it any} pair of maximal ideals in $A$ with the property that 
the finite fields $A/m_1$ and $A/m_2$ have distinct characteristics. Then 
$\mu (M^{4n})$ divides the cardinality of the finite group $GL(2k+1, A/m_1)\times GL(2k+1, A/m_2)$. 

\vskip 5pt

\noindent
{\bf Proof:} Given such a subring and a pair of maximal ideals, Deligne and Sullivan \cite{ds} exhibit a finite cover 
$\bar{M}^{4n}$ of $M^{4n}$ having the property that:
\vskip 2pt
\noindent
(1) the pullback bundle to $\bar{M}^{4n}$ is trivial,\\
(2) the degree of the cover divides $|GL(2k+1, A/m_1)\times GL(2k+1, A/m_1)|$. 
\vskip 2pt
\noindent
But Okun shows (\cite{bo}, proof of theorem 5.1), that there is a tangential map from $\bar{M}^{4n}$ to 
$M_U$, hence applying Corollary 6 completes our proof.

\vskip 5pt

\noindent
{\it Remark:}  The previous corollary tells us that, in some sense, the {\it complexity} of the representation
$\Gamma \rightarrow G_C\subset GL(k,\mathbb C)$ can be estimated from below in terms of the Pontrjagin numbers
of the quotient $\Gamma \backslash G/K$.

\vskip 10pt

\noindent
{\bf 5. Some open questions.}

\vskip 10pt

There remain a few interesting questions along the line of inquiry we are considering.  For starters, Okun has
provided sufficient conditions for establishing non-zero degree of the tangential map he constructs.  One can 
ask the:

\vskip 5pt 
\noindent 
{\bf Question:} Are there examples where Okun's tangential map has zero degree?  In particular, if one has a 
locally symmetric space modelled on $SL(n,\mathbb R)/SO(n)$, does the tangential map to the dual $SU(n)/SO(n)$
have non-zero degree?

\vskip 5pt

Of course, the interest in the special case of $SL(n,\mathbb R)/SO(n)$ is due to the ``universality'' of this
example: every other locally symmetric space of non-positive curvature isometrically embedds in a space modelled
on $SL(n,\mathbb R)/SO(n)$.  Now note that while the relationship between the cohomologies of $M^n$ and $M_U$ (with
real coefficients) is well understood (and has been much studied) since the work of Matsushima \cite{mt}, virtually 
nothing is known about the relationship between the cohomologies with other coefficients.  One can ask:

\vskip 5pt
\noindent 
{\bf Question:} If $t:M^n\rightarrow M_U$ is the tangential map, what can one say about the induced map $t^*:H^*(M_U,
\mathbb Z_p)\rightarrow H^*(M,\mathbb Z_p)$?  

\vskip 5pt

In particular, the case where $p=2$ would be of some particular interest, as the Stiefel-Whitney classes lie in these
cohomology groups.  Finally, we point out that there are other classes of non-positively curved Riemmanian manifolds,
arising from Schroeder's cusp closing construction (\cite{sch}, \cite{hs}), doubling constructions, and related 
techniques.  

\vskip 5pt
\noindent
{\bf Question:}  Compute the characteristic classes for the remaining known examples of non-positively curved manifolds.
\vskip 10pt

\noindent
{\bf Acknowledgments}
\vskip 10pt

The second author would like thank the first for the exchange of ideas that made this collaboration possible.
The authors would also like to thank  Professor F. T. Farrell for
his suggestions during the course of this work, and to thank C. Connell and F. Raymond for some helpful
e-mails.

\end{document}